\documentclass[10pt]{article}
\usepackage{latexsym}
\usepackage{amsfonts}
\usepackage{enumerate}
\usepackage{multicol}
\usepackage{graphicx}
\usepackage{amssymb}
\usepackage{amsmath}
\usepackage{epic}

\title{On two questions about restricted sumsets in finite abelian groups  \\[.4in]}

\author{B\'{e}la Bajnok \\[.1in] {\small Department of Mathematics, Gettysburg College} \\
{\small 300 N. Washington Street, Gettysburg, PA 17325-1486 USA} \\{\small E-mail:  bbajnok@gettysburg.edu} \\ [.2in]
and \\[.2in]
Samuel Edwards \\[.1in] {\small Department of Mathematics, Gettysburg College} \\
{\small 300 N. Washington Street, Gettysburg, PA 17325-1486 USA} \\{\small E-mail:  edwasa01@gettysburg.edu} \\ [.2in]
 \\[.4in]}

\date{July 19, 2016}

\newtheorem{thm}{Theorem}

\newtheorem{lem}[thm]{Lemma}
\newtheorem{cor}[thm]{Corollary}
\newtheorem{prop}[thm]{Proposition}

\begin{document}

\maketitle

\begin{abstract}
 Let $G$ be an abelian group of finite order $n$, and let $h$ be a positive integer.  A subset $A$ of $G$ is called {\em weakly $h$-incomplete}, if not every element of $G$ can be written as the sum of $h$ distinct elements of $A$; in particular, if $A$ does not contain $h$ distinct elements that add to zero, then $A$ is called {\em weakly $h$-zero-sum-free}.  We investigate the maximum size of weakly $h$-incomplete and weakly $h$-zero-sum-free sets in $G$, denoted by $C_h(G)$ and $Z_h(G)$, respectively.  Among our results are the following: (i) If $G$ is of odd order and $(n-1)/2 \leq h \leq n-2$, then $C_h(G)=Z_h(G)=h+1$, unless $G$ is an elementary abelian 3-group and $h=n-3$; (ii)  If $G$ is an elementary abelian 2-group and $n/2 \leq h \leq n-2$, then $C_h(G)=Z_h(G)=h+2$, unless $h=n-4$.  
      
\end{abstract}



\thispagestyle{empty}

\section{Introduction}

Throughout this paper, $G$ denotes a finite abelian group of order $n \geq 2$, written in additive notation.  As is well known, $G$ has a unique {\em invariant decomposition}: that is, it can be written uniquely as the direct product of nontrivial cyclic terms with the order of each term dividing the order of the next; we let $q$ and $r$ denote the {\em exponent} (the order of the last term) and {\em rank} (the number of terms) of $G$, respectively.  If $G$ is {\em cyclic}, we identify it with $\mathbb{Z}_n=\mathbb{Z}/n \mathbb{Z}$; more generally, if $G$ is {\em homocyclic}, we write $G=\mathbb{Z}_q^r$.   We let $L$ denote the subset consisting of the identity element of $G$ as well as of all involutions in $G$: that is, $L$ contains all elements of $G$ of order 1 or 2.  Note that $L$ is a subgroup of $G$; in fact, $L$ is isomorphic to the elementary abelian 2-group whose rank equals the number of even-order terms in the invariant decomposition of $G$.

For a subset $A$ of $G$ we let $|A|$ denote the size of $A$ and $s(A)$ denote the sum of the elements of $A$.  For a positive integer $h$, the (unrestricted) {\em $h$-fold sumset} of $A$, denoted by $hA$, is the collection of all elements of $G$ that can be written as the sum of $h$ (not necessarily distinct) elements of $A$, and the {\em $h$-fold restricted sumset} of $A$, denoted by $h \hat{\;} A$, consists of the elements of $G$ that can be written as the sum of $h$ distinct elements of $A$.

Many questions in additive combinatorics focus on properties of sumsets; for example: How large can a subset of $G$ be without its sumset yielding all of $G$?  While the answer to this question is solved for  unrestricted sumsets (see Theorem \ref{Bajnok thm critical} below), we know much less about restricted sumsets.   The two questions we address in this paper are as follows:
\begin{itemize}
\item How large can a subset $A$ of $G$ be if $h \hat{\;} A \neq G$?
\item How large can a subset $A$ of $G$ be if $0 \not \in h \hat{\;} A$?
\end{itemize} 
In particular, we are interested in finding the quantities
$$C_h(G)=\max \{ |A| \mid A \subseteq G, h \hat{\;} A \neq G\}$$ and
$$Z_h(G)=\max \{ |A| \mid A \subseteq G, 0 \not \in h \hat{\;} A\}.$$  
We say that a subset $A$ of $G$ is {\em weakly $h$-incomplete} if $h \hat{\;} A \neq G$ and that $A$ is {\em weakly $h$-zero-sum-free} if $0 \not \in h \hat{\;} A$.  

These questions can be traced back to the paper \cite{ErdHei:1964a} of Erd\H{o}s and Heilbronn, and variations have been investigated by several authors, including Balandraud \cite{Bal:2012a}; Gao and Geroldinger \cite{GaoGer:2006a}; Lev \cite{Lev:2002a}; Nguyen, Szemer\'edi, and Vu \cite{NguSzeVu:2008a}; and Nguyen and Vu \cite{NguVu:2009a}.   (The terms `$h$-incomplete' and `$h$-zero-sum-free'  have been used in the literature, though we added the word `weakly' to signify the fact that we are considering restricted sumsets.)

One particularly well-researched special case is the problem of finding the largest weakly $3$-zero-sum-free sets in the elementary abelian 3-group $\mathbb{Z}_3^r$, as it corresponds to {\em cap sets} in affine geometry; see  \cite{GaoTha:2004a} by Gao and Thangadurai and its references for $r \leq 5$ and \cite{Pot:2008a} by Potechin for the case $r=6$.  The fact that $Z_3(\mathbb{Z}_3^r)$ is only known for $r \leq 6$ cautions us about the extreme difficulty of these questions; in his blog \cite{Tao:2007a}, Tao writes ``Perhaps my favourite open question is the problem on the maximal size of a cap set.''

At the present time, the only type of group for which $Z_h(G)$ and $C_h(G)$ are known for every value of $h$ is the cyclic group of prime order, and this is due to the fact this is the only case when tight lower bounds for the size of $h$-fold restricted sumsets are known.  Namely, solving a thirty-year open question of Erd\H{o}s and Heilbronn, in 1994 Dias Da Silva and Hamidoune  \cite{DiaHam:1994a} proved that in the cyclic group of prime order $p$, for any nonempty subset $A$ and positive integer $h \leq |A|$ we have
$$|h \hat{\;} A| \geq \min \{p, h|A|-h^2+1\}.$$  ({Soon after, Alon, Nathanson, and Ruzsa  provided a different proof; cf.  \cite{AloNatRuz:1995a} and  \cite{AloNatRuz:1996a}.)   
The fact that this bound is tight can be seen by realizing that equality holds when $A$ is an interval (or, more generally, an arithmetic progression): one can readily verify that if $A$ is an interval of size $m$ in $\mathbb{Z}_p$ (with $m \geq h$), then $h \hat{\;} A$ is an interval of size $\min \{p, hm-h^2+1\}.$  Consequently, in $\mathbb{Z}_p$, the maximum size of a weakly $h$-incomplete set is given by the largest integer $m$ for which $hm-h^2+1$ is less than $p$, or $m=\left \lfloor (p-2)/h \right \rfloor +h$.  Furthermore, for this value of $m$, assuming also that $h<p$, we can choose an interval $A$ in $\mathbb{Z}_p$ of size $m$ for which the interval $h \hat{\;} A$ avoids zero.  Therefore, we have the following:  

\begin{thm}
For any prime $p$ and positive integer $h \leq p-1$ we have $$C_h(\mathbb{Z}_p)=Z_h(\mathbb{Z}_p)=\left \lfloor (p-2)/h \right \rfloor +h.$$

\end{thm}

We make the following observation: When $$(p-1)/2 \leq h \leq p-2,$$ then $\left \lfloor (p-2)/h \right \rfloor=1$, and thus $$C_h(\mathbb{Z}_p)=Z_h(\mathbb{Z}_p)=h+1.$$  One goal of this paper is to prove that the same equations hold in almost every group of odd order.  Namely, we prove the following: If $G$ is a group of odd order $n$ that is not an elementary abelian 3-group, and $h$ is an integer with $$(n-1)/2 \leq h \leq n-2,$$ then $$C_h(G)=Z_h(G)=h+1.$$ More generally:

\begin{thm}  \label{thm summary C and Z}
Let $G$ be an abelian group of order $n$ and exponent $q$, and suppose that its subgroup of involutions $L$ has order $l$.  Then for every integer  
$h$ with 
$$(n+l)/2-1 \leq h \leq n-2,$$  we have $$C_h(G)=Z_h(G)=h+1,$$ with the following two exceptions:
\begin{itemize}
  \item If $h=n-3$ and $q=3$, then $C_h(G)=h+1$ and $Z_h(G)=h$.  
\item If $h=n-2$, $l=2$, and $q \equiv 2$ mod 4, then $C_h(G)=h+1$ and $Z_h(G)=h$. 
\end{itemize}
\end{thm}  

Note that Theorem \ref{thm summary C and Z} is vacuous if (and only if) $G$ is an elementary abelian 2-group; for this case we have the following result:

\begin{thm}  \label{thm 2-group summary C and Z}
Let $G$ be an elementary abelian 2-group of order $n=2^r$, and suppose that $h$ is an integer with 
$$n/2-1 \leq h \leq n-2.$$  Then $$C_h(G)=Z_h(G)=h+2,$$ except when $h=n-4$, in which case $C_h(G)=h+2$ and $Z_h(G)=h$. 

\end{thm}  

Given our theorems above---as well as related results such as those in \cite{NguSzeVu:2008a} by Nguyen, Szemer\'edi, and Vu---we may get the impression that $C_h(G)$ and $Z_h(G)$ are usually equal or that at least they are close to one another.  The following example shows that, actually, $C_h(G)$ and $Z_h(G)$ may be arbitrarily far from one another. 

We say that an $m$-subset $A$ of $G$ is a {\em weak Sidon set} in $G$, if $2 \hat{\;} A$ has size exactly ${m \choose 2}$; in other words, if no element of $G$ can be written as a sum of two distinct elements of $A$ in more than one way (not counting the order of the terms).  Weak Sidon sets were introduced and studied by Ruzsa in \cite{Ruz:1993a}; though the same concept under the name ``well spread set'' was investigated earlier; cf.~\cite{Kot:1972a}.     

\begin{prop}  \label{weakly 4-zero-sum-free iff weak Sidon}
Let $G$ be an elementary abelian 2-group.  Then a subset $A$ of $G$ is weakly 4-zero-sum-free if, and only if, it is a weak Sidon set.

\end{prop}

{\em Proof:}  Let us suppose first that $A$ is weakly 4-zero-sum-free in $G$, and that $a_1+a_2=a_3+a_4$ for some elements $a_1, a_2, a_3,$ and $a_4$ of $A$ with $a_1 \neq a_2$ and $a_3 \neq a_4$.  We then have $$a_1+a_2+a_3+a_4=0,$$ which can only happen if the four terms are not pairwise distinct.  By our assumption, this leads to $\{a_1,a_2\} = \{a_3,a_4\}$, which proves that $A$ is a weak Sidon set in $G$.  The other direction is similar.  $\Box$

According to Proposition \ref{weakly 4-zero-sum-free iff weak Sidon}, if $A$ is a weakly 4-zero-sum-free subset of size $m$ in an elementary abelian 2-group $G$ of order $n=2^r$, then $${m \choose 2} \leq n.$$  On the other hand, we clearly have $C_4(G) \geq n/2$.  This yields the following result:

\begin{prop}
Let $G$ be an elementary abelian 2-group of rank $r$.  We then have $$\lim_{r \rightarrow \infty} \left(C_4(G)-Z_4(G) \right)= \infty.$$

\end{prop}

\section{Weakly $h$-incomplete sets}  \label{section on h-crit}

In this section we study the function $$C_h(G)=\max \{ |A| \mid A \subseteq G, h \hat{\;} A \neq G\};$$ but first, we must mention that the related quantity  $$c_h(G)=\max \{ |A| \mid A \subseteq G, h  A \neq G\}$$ has been completely determined in \cite{Baj:2014a}.  The result can be stated as follows:

\begin{thm} [Bajnok; cf.~\cite{Baj:2014a}]  \label{Bajnok thm critical}

For any abelian group $G$ of order $n$ and for every positive integer $h$, we have
$$c_h(G)=\max \left \{ \left (  \left \lfloor (d-2)/h \right \rfloor +1 \right) \cdot n/d\right \},$$ where the maximum is taken over all divisors $d$ of $n$.

\end{thm}

Observe that---unlike $C_h(G)$---the value of $c_h(G)$ depends only on the order $n$ of $G$ and not on its structure.  

Below, we will employ the fact that $c_h(G)$ is known in the case when $G$ has even order.  Namely, by letting 
$$f_h(d)= \left (  \left \lfloor (d-2)/h \right \rfloor +1 \right) \cdot n/d,$$ we see that $f_h(1)=0$, $f_h(2)=n/2$, and for $d \geq 3$, we get
$$f_h(d) \leq \left ( (d-2)/h  +1 \right) \cdot n/d = \left (  (h-2)/d  +1 \right) \cdot n/h \leq \left (  (h-2)/3  +1 \right) \cdot n/h  \leq n/2.$$  
Therefore, we have the following:

\begin{cor}  \label{c_h for even n}
For any abelian group $G$ of even order $n$ and for every integer $h \geq 2$, we have $c_h(G)=n/2$.
\end{cor}

Let us now turn to the function $C_h(G)$.  These values are easy to find for $h=1$, $h=n-1$, and $h=n$:

\begin{prop}  \label{C for easy h}
For any abelian group $G$ of order $n$ we have $C_1(G)=n-1$, $C_{n-1}(G)=n-1$, and $C_n(G)=n$.

\end{prop}

{\em Proof:}  Each of these claims is quite obvious; for example, to see that $C_{n-1}(G)=n-1$, note that for any subset $A$ of $G$ of size $n-1$, $(n-1) \hat{\;} A$ consists of a single element, and, on the other hand, $(n-1) \hat{\;}G=G$, since for each $g \in G$ we have $s(G \setminus \{s(G)-g\})=g.$  $\Box$

Next, we establish $C_h(G)$ for $h=2$:

\begin{thm}  \label{C_2}
Let $G$ be an abelian group of order $n$, and suppose that its subgroup of involutions has order $l$.   We then have $C_2(G)=(n+l)/2$.

\end{thm}

{\em Proof:}  First, we prove that $C_2(G) \geq (n+l)/2$ by finding a subset $A$ of $G$ with $$|A|=(n+l)/2$$ for which $2 \hat{\;} A \neq G$.  Observe that the elements of $G \setminus L$ are distinct from their inverses, so we have a (possibly empty) subset $K$ of $G \setminus L$ with which $$G=L \cup K \cup (-K),$$ and $L$, $K$, and $-K$ are pairwise disjoint.   Now set $A=L \cup K.$  Clearly, $A$ has the right size; furthermore, it is easy to verify that $0 \not \in 2\hat{\;}A$ and thus $2\hat{\;}A \neq G$.

To prove that $C_2(G) \leq (n+l)/2$, we need to prove that for every subset $A$ of $G$ of size larger than $(n+l)/2$, we have $2\hat{\;}A=G$.  Since this trivially holds when $L=G$, we assume that $L \neq G$.

To continue, we need the following lemma.

{\bf Claim:}  For a given $g \in G$, let $L_g=\{x \in G \mid 2x=g\}$.  If $L_g \neq \emptyset$, then $|L_g| = l.$

{\em Proof of Claim:} Choose an element $x \in L_g$.  Then, for every $y \in L_g$, we have $2(x-y)=0$, and thus $x-y \in L$.  Therefore, $x-L_g \subseteq L$, so  $|x-L_g|=|L_g| \leq l.$  Similarly, $x+L \subseteq L_g$, so  $|x+L|=l \leq |L_g|.$  This proves our claim.  

Now let $m = (n+l)/2+1$. Note that our assumption on $G$ implies that $3 \leq m \leq n$.

Let $A$ be an $m$-subset of $G$, let $g \in G$ be arbitrary, and set $B= g-A$.  Then $|B|=m$, and thus
$$|A \cap B|=|A|+|B|-|A \cup B|  \geq 2m-n  = l+2.$$  
By our claim, we must have an element $a_1 \in A \cap B$ for which $a_1 \not \in L_g$.  Since $a_1 \in A \cap B$, we also have an element $a_2 \in A$ for which $a_1=g-a_2$ and thus $g=a_1+a_2$.  But $a_1 \not \in L_g$, and therefore $a_2 \neq a_1$.  In other words, $g \in 2\hat{\;}A$; since $g$ was arbitrary, we have $G=2\hat{\;}A$, as claimed.
$\Box$

The value of $C_3(G)$ is not known in general and is, in fact, the subject of active interest---see \cite{Baj:2014a}.  Here we present the result for elementary abelian 2-groups:

\begin{thm}  \label{C_3 kappa=2}
If $G$ is the elementary abelian 2-group of order $n=2^r$, then  $C_3(G) =n/2+1.$

\end{thm}

{\em Proof:}  Let $H$ be a subgroup of index 2 in $G$, select an arbitrary element $g \in G \setminus H$, and let $A=H \cup \{g\}$.  Clearly, $g \not \in3 \hat{\;} H$; furthermore, since no two distinct elements of $H$ add to zero, we have $g \not \in 3 \hat{\;} A$.  Therefore, $C_3(G) \geq n/2+1.$

Now let $B$ be a subset of $G$ of size $n/2+2$; we need to show that $3 \hat{\;} B =G$.  (This part of our argument is based on the proof of Theorem 1 in \cite{Lev:2002a}.)  Suppose, indirectly, that this is not so.  Let $g \in G \setminus 3 \hat{\;} B$, and $C=(g+B) \setminus \{0\}$.  Since $|C|=|B|-1=n/2+1$, by Corollary \ref{c_h for even n}, we must have $3C=G$, in particular, $0 \in 3C$.  Therefore, we have elements $c_1, c_2,$ and $c_3$ that add to 0, and thus elements $b_1, b_2,$ and $b_3$ in $B$ for which 
$$(g+b_1)+(g+b_2)+(g+b_3)=0.$$  But $2g=0$ in $G$, so we get $g=b_1+b_2+b_3$.  Since $g \in G \setminus 3 \hat{\;} B$, this can only happen if two of $b_1$, $b_2$, or $b_3$, say $b_1$ and $b_2$, equal each other.  Therefore, $b_1+b_2=0$, so $g=b_3$, and thus $g+b_3=0$.  But this is a contradiction, since $0 \not \in C$.  $\Box$

Regarding the general case, we present an immediate lower bound for $C_h(G)$.  Observe that, if $A$ is any subset of size $h+1$ in $G$, then $h \hat{\;} A$ has size $h+1$ as well.  This yields:

\begin{prop}  \label{C is at least h+1}
For any abelian group $G$ of order $n$ and for every positive integer $h \leq n-2$, we have $C_h(G) \geq h+1$.

\end{prop}

We are now ready to establish our results for $C_h(G)$ for `large' $h$.   
The following lemma will prove useful.

\begin{lem}  \label{C_h lower-upper}
Let $G$ be a finite abelian group, and suppose that $m$ and $h$ are integers for which $$C_{h+1}(G) \leq m \leq C_h(G).$$  Then $C_{m-h}(G)=m$.

\end{lem}

{\em Proof:}  Since $m \leq C_h(G)$, there exists a subset $A$ of $G$ of size $m$ for which $h \hat{\;} A \neq G$.  But $(m-h) \hat{\;} A$ and $h \hat{\;} A$ have the same size, so we must have $(m-h) \hat{\;} A \neq G$ as well, and thus $C_{m-h}(G) \geq m$.

Now let $B$ be any subset of $G$ of size $m+1$; we need to prove that $(m-h) \hat{\;} B= G$.  Since $(m-h) \hat{\;} B$ and $(h+1) \hat{\;} B$ have the same size, we can show that $(h+1) \hat{\;} B=G$ instead.  Since that follows from  $C_{h+1}(G) \leq m$, our proof is complete.  $\Box$

According to the following result, our lower bound of Proposition \ref{C is at least h+1} is actually exact when $h$ is `large':

\begin{thm}  \label{C_h(G)=h+1}
Let $G$ be an abelian group of order $n$, and suppose that its subgroup of involutions has order $l$. Then for every integer $h$ with $$(n+l)/2 -1 \leq h \leq n-2,$$  we have $C_h(G)=h+1$.

\end{thm}

{\em Proof:}  Our claim follows from Proposition \ref{C for easy h}, Theorem \ref{C_2}, and Lemma \ref{C_h lower-upper}, since 
$$C_2(G) = (n+l)/2 \leq h+1 \leq n-1 = C_1(G).$$  $\Box$

We should point out that, when the order of $G$ is odd, then $L=\{0\}$, so we have $C_h(G)=h+1$ for all $(n-1)/2 \leq h \leq n-2$.  More generally, when $L \neq G$, then, since $L$ is a subgroup of $G$, $(n+l)/2$ is at most $3n/4$, so Theorem \ref{C_h(G)=h+1} establishes the function $C_h(G)$ for at least when $h \in [3n/4, n-2]$.  However, Theorem \ref{C_h(G)=h+1} is void when $L=G$; in this case we have the following two results:

\begin{thm}  \label{Z_2^r h+2}
Suppose that $G$ is the elementary abelian 2-group of order $n=2^r$.   
\begin{enumerate}
\item For each integer $h$ with $n/2-1 \leq h \leq n-2,$ we have  $C_h(G) =h+2.$
\item For each integer $h$ with $4 \leq h \leq n/2-2,$ we have  $$n/2 \leq C_h(G) \leq n/2 +h-2.$$
\end{enumerate}
\end{thm}

{\em Proof:}  Our first claim follows from Theorem \ref{C_2}, Theorem \ref{C_3 kappa=2}, and Lemma \ref{C_h lower-upper}, since 
$$C_3(G) = n/2+1 \leq h+2 \leq n = C_2(G).$$ 

The first inequality of the second claim follows from Corollary \ref{c_h for even n}, since $c_h(G) \leq C_h(G)$.  To prove the second inequality, let $A$ be a subset of $G$ of size $n/2+h-1$.  Let us fix a subset $B$ of $A$ of size $h-3$, and let $C = A \setminus B$.  Then $C$ has size $n/2+2$, so $3 \hat{\;} C=G$ by Theorem  \ref{C_3 kappa=2}, and thus $(h-3) \hat{\;} B + 3 \hat{\;} C=G$ as well.  But $(h-3) \hat{\;} B + 3 \hat{\;} C \subseteq h \hat{\;} A$, so $h \hat{\;} A=G$, which proves our claim.  $\Box$

\section{Zero-sum sets of given size}

In this section we develop some results that lay the groundwork for our study of $Z_h(G)$ in Section \ref{section on Z_h(G)}.  We believe these results are of independent interest.

We start with the following easy lemma.

\begin{lem} \label{sum G}
Suppose that $G$ is a finite abelian group with $L$ as the subgroup of involutions; let $|L|=l$.  
\begin{enumerate} \item If $l=2$ with $L=\{0,e\}$, then the sum $s(G)$ of the elements of $G$ equals $e$.  \item If $l \neq 2$, then $s(G)=0$. \end{enumerate}

\end{lem}

{\em Proof:} Recall that $L$ is isomorphic to an elementary abelian 2-group, hence $s(L)=0$, unless $l=2$, in which case $s(L)$ equals the unique element of order 2.  
Our claims follow from the fact that we have $s(G)=s(L)$.  $\Box$ 

We now classify all positive integers $m$ for which one can find $m$ nonzero elements in a given abelian group $G$ that add to $0$.  We separate the cases when $G$ is an elementary abelian 2-group and when it is not.

\begin{thm}  \label{thm on zero-sum sets}
Let $G$ be the elementary abelian 2-group of order $n=2^r$, and let $m$ be a positive integer.  Then $G \setminus \{0\}$ contains a zero-sum subset of size $m$ if, and only if, $3 \leq m \leq n-4$ or $m=n-1$.

\end{thm}

{\em Proof:}  For a given positive integer $k$, let $M(k)$ denote the set of nonnegative integers $m$ for which $\mathbb{Z}_2^k \setminus \{0\}$ contains a zero-sum subset of size $m$.  We start by stating and proving three easy claims about $M(k)$.  

{\bf Claim 1:}  Suppose that $k \geq 2$.  We then have $m \in M(k)$ if, and only if, $2^k-1-m \in M(k)$.

{\em Proof of Claim 1:}  Observe that by Lemma \ref{sum G}, $s(\mathbb{Z}_2^k)=0$, and thus $s(\mathbb{Z}_2^k \setminus \{0\})=0$.  Therefore, for any $A \subseteq \mathbb{Z}_2^k \setminus \{0\}$, we have $$s(A)=s((\mathbb{Z}_2^k \setminus \{0\}) \setminus A),$$ from which our claim follows.

{\bf Claim 2:}  If $m \in M(k)$ for some positive integer $k \geq 2$, then $m \in M(k+1)$.

{\em Proof of Claim 2:}  Clearly, if $A$ is a subset of $\mathbb{Z}_2^k \setminus \{0\}$ of size $m$ with $s(A)=0$, then $B=\{0\} \times A$ is a subset of $\mathbb{Z}_2^{k+1} \setminus \{0\}$ of size $m$ with $s(B)=0$.

{\bf Claim 3:}  Let $k$ and $l$ be integers so that $2 \leq l \leq k$.  If $m \in M(k)$, then $m+2^l \in M(k+1)$.

{\em Proof of Claim 3:}  As in the proof of Claim 2, if $A$ is a subset of $\mathbb{Z}_2^k \setminus \{0\}$ of size $m$ with $s(A)=0$, then $B=\{0\} \times A$ is a subset of $\mathbb{Z}_2^{k+1} \setminus \{0\}$ of size $m$ with $s(B)=0$.    

Let $H$ be a subgroup of order $2^l$ in $\mathbb{Z}_2^k$.  Then $C=\{1\} \times H$ is a subset of $\mathbb{Z}_2^{k+1} \setminus \{0\}$ of size $2^l$ with $s(C)=0$.  Therefore, $B \cup C \subseteq \mathbb{Z}_2^{k+1} \setminus \{0\}$ has size $m+2^l$ and $s(B \cup C)=0$, and thus $m+2^l \in M(k+1)$, as claimed. 

We are now ready to prove Theorem \ref{thm on zero-sum sets}.  Suppose that $G$ has rank $r \geq 2$; we need to prove that $$M(r)=\{0\} \cup \{3,4, \dots, 2^r-4\} \cup \{2^r-1\}.$$  We trivially have $0 \in M(r)$ and $1 \not \in M(r)$.  Furthermore, $2 \not \in M(r)$ follows from the fact that each element of $\mathbb{Z}_2^r$ is its own inverse.  By Claim 1, we then have $2^r-3 \not \in M(r)$, $2^r-2 \not \in M(r)$, and $2^r-1 \in M(r)$.

Assume now that $3 \leq m \leq 2^r-4$; we need to prove that $m \in M(r)$.  Our assumption implies that $r \geq 3$; we verify our claim for $r=3$ and $r=4$, then proceed by induction.

Recall that $2^r-1 \in M(r)$ for each $r \geq 2$; in particular, $3 \in M(2)$ and $7 \in M(3)$.  Therefore, by Claim 2, we have $3 \in M(3)$, $3 \in M(4)$, and $7 \in M(4)$.  Furthermore, $3 \in M(3)$ implies that $4 \in M(3)$ by Claim 1, and thus $4 \in M(4)$ by Claim 2.  By Claim 1, we then also have $\{8,11,12\} \subseteq M(4)$.  This completes the case of $r=3$, and leaves only $m=5,6,9,10$ to be verified for $r=4$; by Claim 1, it suffices to do this for $m=5$ and $m=6$. 

For $i \in \{1,2,3,4\}$, we let $e_i$ denote the element of $\mathbb{Z}_2^4$ with a 1 in the $i$-th position and $0$ everywhere else.  Then the sets $$\{e_1,e_2,e_3,e_4,e_1+e_2+e_3+e_4\}$$ and $$\{e_1,e_2,e_3,e_4,e_1+e_2,e_3+e_4\}$$ show that  $5 \in M(4)$ and $6 \in M(4)$.  This completes our claim for $r=4$. 

Suppose now that $k \geq 4$ and $m \in M(k)$ for every $3 \leq m \leq 2^k-4$; we will show that $m \in M(k+1)$ for every $3 \leq m \leq 2^{k+1}-4$.  For $3 \leq m \leq 2^k-4$, this follows from Claim 2.  Since $k \geq 4$, we have $3 \leq 2^k-7$, so $2^k-7 \in M(k)$, and thus $2^k-3 \in M(k+1)$ by Claim 3; similarly, $2^k-2 \in M(k+1)$ and $2^k-1 \in M(k+1)$.  Therefore, $m \in M(k+1)$ for every $3 \leq k \leq 2^{k}-1$, and thus $m \in M(k+1)$ for every $2^k \leq m \leq 2^{k+1}-4$ as well by Claim 1.  This completes our proof.  $\Box$

\begin{thm}  \label{thm on zer-sum sets of nonzeros}
Let $G$ be an abelian group of order $n$ that is not isomorphic to an elementary abelian 2-group.  Suppose that the subgroup of involutions in $G$ has order $l$, 
and let $m$ be a positive integer.  Then $G \setminus \{0\}$ contains a zero-sum subset of size $m$ if, and only if, one of the following conditions holds:
\begin{itemize}
  \item $2 \leq m \leq n-3$;
  \item $m=n-2$ and $l=2$; or
\item $m=n-1$ and $l \neq 2$.
\end{itemize}

\end{thm}

{\em Proof:}  We may clearly assume that $2 \leq m \leq n-1$.  Let us write $\mathrm{Ord}(G,2)=L \setminus \{0\}$ and $$G=\{0\} \cup \mathrm{Ord}(G,2) \cup K \cup -K,$$ where the four components are pairwise disjoint and, since $G$ is not isomorphic to an elementary abelian 2-group, $K$ and $-K$ are nonempty.  We examine three cases.

{\bf Case 1:} $l=1$.

In this case, $q$ and $n$ are odd, and $\mathrm{Ord}(G,2)=\emptyset$, and thus $G \setminus \{0\}=K \cup -K$.  Clearly, $G \setminus \{0\}$ clearly contains a zero-sum subset of every even size $m \leq n-1$.  Furthermore, we see that $G \setminus \{0\}$ does not have a zero-sum set of size $n-2$.  It remains to be shown that $G \setminus \{0\}$ contains a zero-sum subset of every odd size $3 \leq m \leq n-4$.

If $n=7$, then the set $\{1,2,4\}$ proves our claim, so let us assume that $n \geq 9$ or, equivalently, that $|K|\geq 4$.  Let $g_1$ be any element of $K$; since $|K|\geq 4$, we can find another element $g_2 \in K$ so that $g_2 \neq -2g_1$ and $g_2 \neq \frac{q -1}{2} g_1$.

We first prove that the six elements $\pm g_1, \pm g_2$, and $\pm (g_1+g_2)$ are pairwise distinct.  Indeed, $g_1$ and $g_2$ are distinct elements of $K$, so $-g_1$ and $-g_2$ are distinct elements of $-K$.  So $g_1+g_2 \neq 0$, and thus one of $g_1+g_2$ or $-(g_1+g_2)$ is an element of $K$ and the other is an element of $-K$.  If $g_1+g_2$ is in $K$, then it must be distinct from both $g_1$ and $g_2$, since neither of these is $0$, and so $-(g_1+g_2)$ is distinct from $-g_1$ and $-g_2$ as well.   Furthermore, if $g_1+g_2$ is in $-K$, then it must be distinct from $-g_1$ since $g_2 \neq -2g_1$, and  if it were equal to $-g_2$, then we would get $2g_2=-g_1$, so $\frac{q+1}{2} \cdot 2g_2=\frac{q+1}{2} \cdot (-g_1)$, that is, $g_2=\frac{q -1}{2} g_1$, which we ruled out.

Therefore, we are able to partition $G$ as
$$G=\{0\} \cup \{\pm g_1, \pm g_2, \pm (g_1+g_2)\} \cup K' \cup -K',$$ where $K' \subset K$ and $|K'|=(n-7)/2$.  Note that $(m-3)/2 \leq |K'|$;   let $K_1 \subseteq K'$ of  size $(m-3)/2$.  Then $$A=\{g_1, g_2, -(g_1+g_2)\} \cup K_1 \cup -K_1$$ has size $m$ and its elements sum to $0$. 

{\bf Case 2:}  $l=2$.

In this case, $q$ is even and $n/q$ is odd, and $|\mathrm{Ord}(G,2)|=1$.  Let $\mathrm{Ord}(G,2)=\{e\}$; we then have
$$G=\{0\} \cup \{e\} \cup K \cup -K.$$  Clearly, $G \setminus \{0\}$ contains a zero-sum subset of every even size $m \leq n-2$; we consider odd values of $m$ next.

The case  of $m=n-1$ is settled by the fact that the elements of $G \setminus \{0\}$ add up to $e$ by Lemma \ref{sum G}.  Next, we consider $m=n-3$, in which case we are looking for a set $A$ of the form $$A=G \setminus \{0,g_1,g_2\}$$ whose elements add to $0$.  Now $m \geq 3$, so $n \geq 6$, and since $q$ is even and $n/q$ is odd, we then must have $q \geq 6$ as well.  Let $g_1$ be any element of $G$ of order $q$, and let $g_2=e-g_1$.  Then $g_1$ and $g_2$ are distinct nonzero elements of $G$, since $g_1=g_2$ would imply that $g_1$ has order at most 4.  Thus $A$ satisfies our requirements.

This leaves us with the cases of odd $m$ values with $3 \leq m \leq n-5$.  If $n=8$, then our assumptions imply that $G$ is cyclic, in which case the set $\{1,3,4\}$ satisfies our claim.  If $n \geq 10$, then $|K| \geq 4$, so this case can be handled as in Case 1 above.    

{\bf Case 3:}  $l >2$.  

In this case, $q$ and $n/q$ are even, and $|\mathrm{Ord}(G,2)|>1$.  Note that the elements of $G$, and thus the elements of $G \setminus \{0\}$, sum to $0$; this settles the cases of $m=n-1$ and $m=n-2$.  We need to show that a zero-sum subset of $G \setminus \{0\}$ of size $m$ exists for every $2 \leq m \leq n-3$.  

Recall that $L$ is isomorphic to an elementary abelian 2-group, so $|\mathrm{Ord}(G,2)|$ is 1 less than a power of 2; so, by assumption, it equals 3 or is at least 7.

Suppose first that $|\mathrm{Ord}(G,2)|=3$.  Since the three elements of $\mathrm{Ord}(G,2)$ add to $0$, $G \setminus \{0\}$ contains a zero-sum subset of every odd size $3 \leq m \leq 3+2|K|=n-1$.  Clearly, $G \setminus \{0\}$ also contains a zero-sum subset of every even size $3 \leq m \leq 2|K|=n-4$ as well, completing this case.

Suppose now that $|\mathrm{Ord}(G,2)| \geq 7$.  By Theorem \ref{thm on zero-sum sets}, $\mathrm{Ord}(G,2)$, and thus $G \setminus \{0\}$, contains a zero-sum subset of size $m$ for every $2 \leq m \leq |\mathrm{Ord}(G,2)|-3$.  If $|\mathrm{Ord}(G,2)|-2 \leq m \leq n-4$, then we may write $m$ as $m=m_1+2k_1$, with $0 \leq k_1 \leq |K|$ and $m_1=|\mathrm{Ord}(G,2)|-3$ (if $m$ is even) or $m_1=|\mathrm{Ord}(G,2)|-4$ (if $m$ is odd).  Therefore, $G \setminus \{0\}$ contains a zero-sum subset of every size $m$ with $2 \leq m \leq n-4$.  Finally, if $m=n-3$, then $m=|\mathrm{Ord}(G,2)|+2(|K|-1)$, so again $G \setminus \{0\}$ contains a zero-sum subset of size $m$.   This completes our proof.  $\Box$

\begin{cor}  \label{cor on zer-sum sets}
Let $G$ be an abelian group of order $n$.  Suppose that the subgroup of involutions in $G$ has order $l$, and let $m$ be a positive integer with $m \leq n$.  Then $G$ contains a zero-sum subset of size $m$ with the following exceptions:
\begin{itemize}
  \item $G$ is isomorphic to an elementary abelian 2-group and $m \in \{2,n-2\}$; or
  \item $l=2$ and $m=n$.
\end{itemize}

\end{cor}

{\em Proof:}  The claim is trivial for $m=1$, and is a restatement of Lemma \ref{sum G} if $m=n$.  If $G$ and $m$ are such that $G \setminus \{0\}$ contains a zero-sum set $A$ of size $m$ or $m-1$, then either $A$ or $A \cup \{0\}$ satisfies our claim.  By Theorems \ref{thm on zero-sum sets} and \ref{thm on zer-sum sets of nonzeros}, this leaves only the case when $G$ is isomorphic to an elementary abelian 2-group and $m=2$ or $m=n-2$, for which the claim follows from the fact that each element is its own inverse then.  $\Box$

\begin{cor}  \label{cor on sets not containing their sum}
Let $G$ be an abelian group of order $n$.  Suppose that the subgroup of involutions in $G$ has order $l$, and let $m$ be a positive integer.  Then $G$ contains a subset $A$ of size $m$ for which $s(A) \not \in A$ if, and only if, one of the following conditions holds:
\begin{itemize}
  \item $2 \leq m \leq n-4$;
  \item $m=n-3$ and $G$ is not isomorphic to an elementary abelian 2-group;
  \item $m=n-2$ and $G$ is not isomorphic to an elementary abelian 3-group; or
  \item $m=n-1$ and $l \neq 2$; or $m=n-1$, $l=2$, and $q$ is divisible by 4.
\end{itemize}

\end{cor}

{\em Proof:}  We can clearly assume that $2 \leq m \leq n-1$, and by Theorems \ref{thm on zero-sum sets} and \ref{thm on zer-sum sets of nonzeros}, it suffices to consider the following cases:
\begin{enumerate}[(i)]
  \item $m=n-3$ and $G$ is isomorphic to an elementary abelian 2-group;
  \item $m=n-2$ and $l \neq 2$; and
  \item $m=n-1$, $l=2$.  
\end{enumerate}

If $m=n-3$ and $G$ is isomorphic to an elementary abelian 2-group, then an $m$-set $A$ with $s(A) \not \in A$ exists if, and only if, we can find distinct elements $a_1,a_2,$ and $a_3$ in $G$ for which $a_1+a_2+a_3 \in \{a_1,a_2,a_3\}$.  This is not possible, since two distinct elements do not add to $0$.  

The cases to be considered for $m=n-2$ are exactly those where, by Lemma \ref{sum G}, $s(G)=0$.  Therefore, an $m$-set $A$ with $s(A) \not \in A$ exists if, and only if, we can find distinct elements $a_1$ and $a_2$ in $G$ for which $-a_1-a_2 \in \{a_1,a_2\}$, that is, $a_2 \neq -2a_1$ or $a_1 \neq -2a_2$.  This is possible exactly when $G$ has an element whose order is not a divisor of 3.

Finally, suppose that $m=n-1$ and $l=2$.  In this case, by Lemma \ref{sum G}, $s(G)=e$ where $e$ is the unique element of $G$ of order 2. Therefore, an $m$-set $A$ with $s(A) \not \in A$ exists if, and only if, $G$ contains an element $a$ for which $2a=e$, which is possible exactly when $q$ is divisible by 4.  $\Box$

\section{Weakly $h$-zero-sum-free sets}  \label{section on Z_h(G)}

We start by determining $$Z_h(G)=\max \{ |A| \mid A \subseteq G, 0 \not \in h \hat{\;} A\}$$ for $h=1, 2, n-1$, and $n$.

\begin{prop}
Let $G$ be an abelian group of order $n$, and suppose that its subgroup of involutions has order $l$.  We have 
\begin{enumerate} 
  \item $Z_1(G)=n-1$;
  \item $Z_2(G)=(n+l)/2$;
  \item $Z_{n-1}(G)=n-1$;
  \item $Z_n(G)=n$ when $l=2$, and $Z_n(G)=n-1$ when $l \neq 2$.
\end{enumerate}

\end{prop}

{\em Proof:}  The first claim is trivial, since $G \setminus \{0\}$ is weakly 1-zero-sum-free.  Let us write $G=L \cup K \cup (-K).$  Clearly, $A=L \cup K$ is weakly 2-zero-sum-free.  On the other hand, if $B$ has size more than $(n+l)/2$, then it contains at least $(n-l)/2+1=|K|+1$ elements of $K \cup (-K)$, so it is not weakly 2-zero-sum-free.

To prove that $Z_{n-1}(G)=n-1$, let $g=s(G)$.  Then $s(G \setminus \{g\})=0$, so $Z_{n-1}(G) \leq n-1$.  But for every element $g' \in G \setminus \{g\}$, we have $s(G \setminus \{g'\})=g-g' \neq 0$, so $Z_{n-1}(G) \geq n-1$.  Our last claim follows from Lemma \ref{sum G}.  $\Box$

We can easily establish the following lower bound:

\begin{prop}  \label{Z_h >=h}
For any abelian group $G$ of order $n$ and all positive integers $h \leq n-1$ we have $Z_h(G) \geq h$.

\end{prop}

{\em Proof:}  Let $A$ be any subset of $G$ of size $h$.  If $s(A) \neq 0$, we are done.  Otherwise, choose elements $a \in A$ and $b \in G \setminus A$.  Then for $B=(A \setminus  \{a\}) \cup \{b\}$ we have $|B|=h$ and $$s(B)=s(A)-a+b=b-a \neq 0.$$  $\Box$  

Next, we present a necessary and sufficient condition for $Z_h(G)$ to be at least $h+1$:

\begin{prop}  \label{when Z=h+1}
Let $G$ be a finite abelian group and $h$ be a positive integer with $h \leq n-1$.  Then $Z_{h}(G) \geq h+1$ if, and only if, there exists a subset $A$ in $G$ of size $h+1$ for which $s(A) \not \in A$.

\end{prop}

{\em Proof:}  Suppose first that $A$ is a subset of $G$ of size $h+1$ for which $s(A) \not \in A$; we prove that $A$ is weakly $h$-zero-sum-free in $G$.  Let $B$ be any subset of size $h$ of $A$, and let $a$ be the element of $A$ for which $B=A \setminus \{a\}$.  Then $s(B)=s(A)-a$; since $s(A) \not \in A$, we have $s(B) \neq 0$, as claimed.  Therefore, $Z_{h}(G) \geq h+1$.

Conversely, assume that all subsets of $G$ of size $h+1$ contain their sum as an element.  Let $A$ be any subset of $G$ of size $h+1$.  By assumption, $s(A) \in A$; let $B = A \setminus \{s(A)\}$.  Then $B$ has size $h$ and $s(B)=0$, so $A$ is not weakly $h$-zero-sum-free in $G$.  Therefore, $Z_{h}(G) \leq h$.  $\Box$

Our next two results establish the value of $Z_h(G)$ for all `large' $h$.  First, we consider groups with exponent at least three:

\begin{thm}
Let $G$ be an abelian group of order $n$ that is not isomorphic to an elementary abelian 2-group, and suppose that its subgroup of involutions has order $l$.  For every integer  $h$ with $$(n+l)/2-1 \leq h \leq n-2,$$ we have 
$$Z_h(G)=\left\{
\begin{array}{cl}
h & \mbox{if}  \; h=n-3 \; \mbox{and} \; q=3; \; \mbox{or} \\ \\
& h=n-2, \; l=2, \; \mbox{and} \; q \equiv 2 \; \mbox{mod 4}; \\ \\
h+1 & \mbox{otherwise}.  
\end{array}\right.$$

\end{thm}

{\em Proof:}  By Proposition \ref{Z_h >=h} and Theorem \ref{C_h(G)=h+1}, we have $$h \leq Z_h(G) \leq h+1.$$  Thus our claim follows from Proposition \ref{when Z=h+1} and Corollary \ref{cor on sets not containing their sum}.  $\Box$

For groups of exponent two, we have the following result:

\begin{thm}
Suppose that $G$ is isomorphic to an elementary abelian 2-group and has order $n=2^r$, and let $h$ be an integer with $n/2-1 \leq h \leq n-2$.  We then have 
$$Z_h(G)=\left\{
\begin{array}{cl}
h & \mbox{if}  \; h=n-4; \\ \\
h+2 & \mbox{otherwise}. 
\end{array}\right.$$

\end{thm}

{\em Proof:}  By Proposition \ref{Z_h >=h} and Theorem \ref{Z_2^r h+2}, we have $$h \leq Z_h(G) \leq h+2.$$  Therefore, our result will follow from the following two claims.

{\bf Claim 1:}  If $h$ is a positive integer with $h \leq n-2$ and $h \neq h-4$, then $Z_h(G) \geq h+2$.

{\em Proof of Claim 1:}  Let $m=h+2$; we then have $3 \leq m \leq n$ with $m \neq n-2$.  Thus, by Corollary \ref{cor on zer-sum sets}, $G$ contains an $m$-subset $A$ with $s(A)=0$; we will prove that $A$ is weakly $h$-zero-sum-free in $G$.  Let $B$ be any $h$-subset of $A$; we assume that $B=A \setminus \{a_1,a_2\}$.  Since $a_1$ and $a_2$ are distinct, we have $a_1+a_2 \neq 0$, and therefore $$s(B)=s(A)-(a_1+a_2)=a_1+a_2 \neq 0.$$  This proves our claim.

{\bf Claim 2:}  We have $Z_{n-4}(G) \leq n-4$.

{\em Proof of Claim 2:}  Suppose that $A$ is an arbitrary subset of $G$ with $|A|=n-3$; we let $A=G \setminus \{a_1,a_2,a_3\}$.  Note that $a_1, a_2,$ and $a_3$ are pairwise distinct, so no two of them add to zero, and thus $a_1+a_2+a_3 \in A$.  Let $B =A \setminus \{a_1+a_2+a_3\}$.  We then have $$s(B)=s(A)-(a_1+a_2+a_3),$$ where $$s(A)=s(G)-(a_1+a_2+a_3)=a_1+a_2+a_3.$$  Thus $s(B)=0$, which proves our claim.  $\Box$

\end{document}